\documentclass[12pt]{elsarticle}
\usepackage[utf8]{inputenc}
\usepackage[english]{babel}
\usepackage{amsmath}
\usepackage{amssymb}
\usepackage{amsfonts}
\usepackage{mathtools}
\usepackage{dsfont}
\usepackage{rotating}
\usepackage{graphicx}
\usepackage{floatflt,epsfig}
\usepackage{lineno,hyperref}
\usepackage{enumerate}
\usepackage{colortbl}
\usepackage{array,tabularx,tabulary,booktabs}
\usepackage{longtable}
\usepackage{multirow}
\usepackage{wrapfig}
\usepackage{subcaption}
\usepackage{xurl}
\usepackage[table]{xcolor}
\newcolumntype{^}{>{\currentrowstyle}}

\usepackage{amsthm}
\usepackage[left=1in,right=1in,top=1in,bottom=1in]{geometry}

\usepackage{amssymb,amsthm}

\newtheorem{lemma}{Lemma}[section]
\newtheorem{theorem}[lemma]{Theorem}

\newtheorem{remark}[lemma]{Remark}
\newtheorem{example}[lemma]{Example}
\newtheorem{construction}[lemma]{Construction}

\def\Sz{{\operatorname{Sz}}}
\def\Ree{{\operatorname{Ree}}}

\journal{Designs, Codes and Cryptography}
\begin{document}
\renewcommand{\abstractname}{Abstract}
\renewcommand{\refname}{References}
\renewcommand{\tablename}{Table}
\renewcommand{\arraystretch}{0.9}
\thispagestyle{empty}
\sloppy

\begin{frontmatter}
\title{Divisible design graphs from a construction of antipodal distance-regular graphs based on Ree groups and Suzuki groups}

\author[02]{Bart De Bruyn}
\ead{Bart.DeBruyn@UGent.be}

\author[01]{Sergey Goryainov}
\ead{sergey.goryainov3@gmail.com}

\author[03]{Yilin Xie}
\ead{xyl_crazystu@hebtu.edu.cn}

\address[01] {School of Mathematical Sciences, Hebei International Joint Research Center for Mathematics and Interdisciplinary Science, Hebei Key Laboratory of Computational Mathematics and Applications, Hebei Workstation for Foreign Academicians, Hebei Normal University, Shijiazhuang  050024, P.R. China}
\address[02] {Ghent University, Department of Mathematics, Computer Science and Statistics, Krijgslaan 299 - S9, 9000 Gent, Belgium}
\address[03] {School of Mathematical Sciences, Hebei Normal University, Shijiazhuang 050024, P.R. China}


\begin{abstract}
In this paper, we present a new construction of divisible design graphs with new parameters, obtained by plugging a difference set into a known construction of antipodal distance-regular graphs of diameter 3 that is based on Ree groups and Suzuki groups.
\end{abstract}

\begin{keyword}
divisible design graph; distance-regular graph; difference set; Ree group; Suzuki group 
\vspace{\baselineskip}
\MSC[2020] 20D06\sep 05B05\sep 05B10\sep 05C25 \sep 05E30 
\end{keyword}
\end{frontmatter}

\section{Introduction} 
A \emph{divisible design graph} (a DDG for short) with parameters $(v,k,\lambda_1,\lambda_2,m,n)$ is a $k$-regular graph on $v$ vertices such that its vertex set can be partitioned into $m$ classes of size $n$ with the following two properties: any two distinct vertices from the same class have precisely $\lambda_1$ common neighbours and any two vertices from different classes have precisely $\lambda_2$ common neighbours. The partition from the definition of a DDG is called the \emph{canonical partition}. DDGs were introduced in \cite{HKM11,M08} as a bridge between graph theory and design theory (the adjacency matrix of a divisible design graph can be regarded as the incidence matrix of a group divisible design \cite{B77}) and have been studied in \cite{BDG25,BG24,CH14,CS22,DGHS24,GK24,GK25,GHKS19,GHKL25,K22,K23,K26,KS21,MR26,P22,PS22,R26,S21,T22}.
Similarly as strongly regular graphs (SRGs), DDGs form an important subclass of Deza graphs \cite{E99}. 

In \cite{DGMX26}, De Bruyn, Goryainov, Ma and Xie construct new families of DDGs which are Cayley graphs of elementary abelian $2$-groups by plugging a difference set of a quotient group into a known construction of antipodal distance-regular graphs of diameter $3$. 
Inspired by this idea, we present in Construction~\ref{constr} and Theorem~\ref{thm:DDG} a construction of new DDGs arising from coset graphs of Suzuki and Ree groups, by plugging a cyclic difference set into Tsiovkina's construction of antipodal distance-regular graphs of diameter $3$ (see \cite{T15}). 
In particular, by plugging Singer difference sets and Paley difference sets into Tsiovkina's construction, we obtain four infinite families of DDGs as examples of our construction.
Several DDGs in the literature are Cayley graphs (see \cite{CS22,DGMX26,KS21,R26}), and Cayley graphs form a subclass of coset graphs. 
To our knowledge, our examples are the first explicit constructions of DDGs arising from coset graphs of simple groups.

This paper is organized as follows. In Section \ref{sec:prelim}, we give preliminary definitions and results on distance-regular graphs, coset graphs, Suzuki and Ree groups, difference sets and some basic number-theoretic results. In Section \ref{sec:construction}, we present our new construction, see Construction~\ref{constr} and Theorem~\ref{thm:DDG}, and conclude with examples. 

\section{Preliminaries}\label{sec:prelim}
In this section we give preliminary definitions and results.

\subsection{Distance-regular graphs}
A \emph{distance-regular graph} is a connected regular graph such that for any two vertices $x$ and $y$, the number of vertices at distance $j$ from $x$ and at distance $\ell$ from $y$ depends only upon $j$, $\ell$, and the distance between $x$ and $y$.
The \emph{intersection array} of a distance-regular graph is the array $(b_{0},b_{1},\ldots,b_{d-1};c_{1},\ldots ,c_{d})$ in which $d$ is the diameter of the graph and for each $0\leq j\leq d$, $b_{j}$ gives the number of neighbours of $y$  at distance $j + 1$ from $x$, $c_{j}$ gives the number of neighbours of $y$ at distance $j-1$ from $x$ for any pair of vertices $x$ and $y$ at distance $j$. There is also the number $a_j$ that gives the number of neighbours of $y$ at distance $j$ from $x$. The numbers $a_j$, $b_j$, $c_j$ are called the \emph{intersection numbers} of the distance-regular graph. They satisfy the equation $a_{j}+b_{j}+c_{j}=k$, where $k = b_0$ is the valency, that is, the number of neighbours, of any vertex.

A distance-regular graph of diameter $d$ is called \emph{antipodal} if the relation on its vertex set defined by the rule ``to coincide or to be at distance $d$'' is an equivalence relation; the equivalence classes are then called \emph{antipodal classes}. Note that whether a distance-regular graph is antipodal or not can be told from its intersection array (see \cite[Proposition 4.2.2(ii)]{BCN89}).

\subsection{Coset graphs}\label{sec:cosgraph}
Here we introduce the notion of coset graphs associated with groups, which will be used later.
Let $G$ be a finite group and $H$ be a subgroup of $G$. 
Let $X_0\subseteq G$ be a subset, and $X=\bigcup_{a\in X_0}HaH$ be the union of $(H,H)$-double cosets represented by $X_0$.
Assume further that $X=X^{-1}$ and $1\notin X$ (the latter is equivalent to $X_0 \cap H= \emptyset$).
Then $\Gamma(G,H,X)$ is an undirected, loopless graph whose vertex set is the set of right cosets of $H$ in $G$, with adjacency defined by 
\[
Hg_1\sim Hg_2
\quad\Longleftrightarrow\quad
g_1g_2^{-1}\in X.
\] 
In particular, if $H=1$, then $X=X_0$, and $\Gamma=\operatorname{Cay}(G,X_0)$ is a Cayley graph.

The following lemma is well known and straightforward to prove, so we omit its proof.
\begin{lemma}\label{lem:trans}
    Let $\Gamma=\Gamma(G,H,X)$.
    For every $g\in G$, the map 
    $\iota_g:V(\Gamma)\to V(\Gamma)$ defined by $\iota_g(Hx)=Hxg$ for each $x\in G$ is an automorphism of $\Gamma$.
\end{lemma}
We call $\iota_g$ the action of $g$ on $V(\Gamma)$ by right multiplication. 
Clearly, $G$ acts transitively on $V(\Gamma)$ by right multiplication.

\subsection{Basic facts about Suzuki and Ree groups}\label{subsec:groups}
Suzuki and Ree groups are two families of finite simple groups of Lie type.
We will construct divisible design graphs as coset graphs of these groups.

Here we list some useful facts about these groups that will be used later. 

\paragraph{Standing assumptions}
Throughout this paper, let $G$ be either $\Sz(q)$ or $\Ree(q)$, where 
\[q=p^{2n+1}>3\] 
for some positive integer $n$. Here $p$ is the defining characteristic, that is, $p=2$ if $G=\Sz(q)$ and $p=3$ if $G=\Ree(q)$. 
Then $|G|=q^2(q^2+1)(q-1)$ if $G=\Sz(q)$ (see for example \cite[Theorem 7.3.3]{JHRLMS}); and $|G|=q^3(q^3+1)(q-1)$ if $G=\Ree(q)$ (see for example \cite{LN85}).
We take $S$ to be a Sylow $p$-subgroup of $G$.
Note that $G$ is simple, so the subgroup generated by $t^G$ for any involution $t$ of $G$ is $G$ itself. 
Suppose that $N_G(S)$ contains all involutions in $G$. 
Then $N_G(S)$ would contain the subgroup generated by $t^G$, forcing $N_G(S)=G$, which is absurd.
It follows that there exists an involution $g\in G \setminus N_G(S)$.
We fix an involution \[g\in G\setminus N_G(S).\] Let $H$ be a subgroup of $N_G(S)$ of odd index $r>1$, where $r$ divides $q-1$. 

\begin{lemma}\label{lem:sylownormalizer}
Under the standing assumptions above, there exists \(f\in N_G(S)\) of order \(q-1\), and the following statements hold
\begin{enumerate}[\normalfont (i)]
    \item if $G=\Sz(q)$, then $|S|=q^2$; if $G=\Ree(q)$, then $|S|=q^3$.
    \item $N_G(S)=S \rtimes \langle f \rangle$.
    \item $|G : N_G(S)|=|S|+1$.
    \item $N_G(S) \cap (N_G(S))^g=\langle f \rangle$, and $f^g=f^{-1}$.
    \item $H=S \rtimes \langle f^r \rangle$, and $H$ is the unique subgroup of index $r$ in $N_G(S)$.
    \item $\{1,f,\cdots, f^{r-1}\}\cup \{f^igs\mid 0\leq i \leq r-1, s\in S\}$  
    forms a complete set of representatives of the right cosets of $H$ in $G$.
    \item for each integer $i$, we have $Hf^igS=Hf^igH$, and \[\{ 1,f,f^2,\cdots, f^{r-1} \}\cup \{g, fg,f^2g,\cdots, f^{r-1}g\}\] forms a complete set of representatives of the double cosets $H \backslash G /H$.
\end{enumerate}
\end{lemma}
\begin{proof}
    \rm(i) follows directly from the order of $G$. \rm(ii)-\rm(iv) can be deduced from \cite[Propositions 3 and 4]{T15}. (In \cite{T15}, they denote $N_G(S)$ as $M$, and $\langle f \rangle$ as $K$.)
    In particular, \rm(ii) follows from Proposition 3 (1) and Proposition 4 (1); \rm(iii) is a direct consequence of Proposition 3 (2) and Proposition 4 (2); and \rm(iv) is a combination of Proposition 3 (1),(3) and Proposition 4 (1),(3). 
    It remains to prove \rm(v), \rm(vi) and \rm(vii).
   
    For \rm(v), we first show the uniqueness of $H$. To this end, let $H_1$ be a subgroup of index $r$ in $N_G(S)$; we proceed to prove that $H_1=H$. 
    Indeed, since $S$ is the unique normal Sylow $p$-subgroup of $N_G(S)$ and $|S|$ divides $|H_1|$, we have $S \mathrel{\unlhd} H_1$.
    Similarly, $S\unlhd H$.
    On the other hand, $H_1/S$ and $H/S$ are both index $r$ subgroups of the cyclic group $N_G(S)/S$, which has a unique subgroup of index $r$.
    It follows that $H_1/S=H/S$, and hence $H_1=H$. 
    Furthermore, it is clear that $S\rtimes \langle f^r\rangle$ is a subgroup of index $r$ of $N_G(S)$, yielding that $H=S\rtimes \langle f^r \rangle$.

    For \rm(vi), 
    \begin{equation}\label{eqngs}
    G= N_G(S) \sqcup \bigsqcup_{s\in S}  N_G(S)gs    
    \end{equation}
    by \cite[Proposition 3(2) and Proposition 4(2)]{T15}.
    On the other hand, we have $N_G(S)=\bigsqcup_{j=0}^{r-1}Hf^j$.
    It follows that \[G=\bigsqcup_{i=0}^{r-1}Hf^i \sqcup (\bigsqcup_{i=0}^{r-1}\bigsqcup_{s\in S}Hf^igs),\] as required.
    
    For \rm(vii), we first show that for each integer $i$, we have $Hf^igH=Hf^igS$.
    To see this, note that $H=S\langle f^r \rangle=\langle f^r \rangle S $, which is due to $S\mathrel{\unlhd}N_G(S)$.
    Hence,
    \[
  \begin{aligned}
   Hf^igH
   &=Hf^igS\langle f^r\rangle
   =Hf^ig\langle f^r\rangle S \\
   &=Hf^i\langle f^{-r}\rangle gS
   =H\langle f^{-r}\rangle f^igS
  =Hf^igS.
   \end{aligned}
 \]
    where the third to last equality holds from the fact $f^g=f^{-1}$ in \rm(iv).
    Next, we compute the double cosets decomposition $H\backslash G/H$. 
    Note that $Hf^jH=HHf^j=Hf^j$, which is because $H\mathrel{\unlhd} N_G(S)$.
    Applying \rm(vi), we have
    \[ G=\bigsqcup _{j=0}^{r-1}Hf^j \sqcup \bigsqcup_{j=0}^{r-1}Hf^jgS = \bigsqcup _{j=0}^{r-1}Hf^jH \sqcup \bigsqcup_{j=0}^{r-1}Hf^jgH ,\] as required.

    This completes the proof.
\end{proof}

The following theorem gives an infinite family of antipodal distance-regular graphs of diameter 3. We remark here that ${}^2G_2(q)$ is the Lie notation for $\Ree(q)$, which is used in Theorem~\ref{thm:DRG1}; while the Lie notation for Suzuki group is ${}^2B_2(q)$. (The notation ${}^2B_2(q)$ is not used in Theorem~\ref{thm:DRG1}.)
Recall that we are assuming $r>1$; the case $r=1$ would give a complete graph
and is therefore excluded.

\begin{theorem}[{\cite[Theorem 1]{T15}}]\label{thm:DRG1}
Let $G \in \{\Sz(q), {}^2G_2(q)\}$, where $q = p^{2n+1} > 3$, $S \in \operatorname{Syl}_p(G)$ with
$p = 2$ if $G = \Sz(q)$, and $p = 3$ if $G = {}^2G_2(q)$. Let $g$ be an involution of $G$ not
contained in $N_G(S)$, let $\langle h\rangle$ be a subgroup of $N_G(S) \cap N_G(S)^g$ of odd index $r > 1$
and $H = S\langle h\rangle$. Then $\Gamma(G, H, HgH)$ is an arc-transitive antipodal distance-regular
graph with intersection array $\{q^2, (q^2-1)(r-1)/r, 1; 1, (q^2-1)/r, q^2\}$ if $G = \Sz(q)$,
or $\{q^3, (q^3 - 1)(r - 1)/r, 1; 1, (q^3 - 1)/r, q^3\}$ if $G = {}^2G_2(q)$, and $\Gamma(G, H, HgH)$
does not depend (up to isomorphism) on the choice of the involution $g \in G \setminus N_G(S)$.   
\end{theorem}

\subsection{Difference sets}
A $(v,k,\lambda)$ \emph{difference set} is a subset $D$ of size $k$ in a group $\mathcal{G}$ of order $v$ such that every non-identity element of $\mathcal{G}$ can be expressed as a product $d_1\ast d_2^{-1}$ of elements of $D$ in exactly $\lambda$ ways, that is, when
$$
D\ast D^{-1} = k\cdot\{e\} + \lambda\cdot (\mathcal{G} \setminus \{e\})
$$
holds.
Note that for a $(v,k,\lambda)$ difference set $D$ in a group $\mathcal{G}$, the complementary set $\mathcal{G} \setminus D$ is a $(v,v-k, v-2k+\lambda)$ difference set.
For any group $\mathcal{G}$ of order $v$ and any element $a \in \mathcal{G}$,
the set $\{a\}$ is a $(v,1,0)$ difference set and the set $\mathcal{G} \setminus \{a\}$ is a $(v,v-1,v-2)$ difference set;
such difference sets are called \emph{trivial}. 

\subsubsection{Singer difference sets}
Let $q$ be a prime power and $m$ be an integer with $m>1$. Let $\mathbb{F}_{q^m}$ (resp. $\mathbb{F}_q$) be the finite field defined on $q^m$  (resp. $q$) elements. Then $\mathbb{F}_q\subseteq\mathbb{F}_{q^m}$.
Define $\mathcal{G}$ to be the quotient group $\mathbb{F}^\times_{q^m}/\mathbb{F}^\times_q$, where $\mathbb{F}^\times_{q^m}=\mathbb{F}_{q^m}\setminus\{0\}$ and $\mathbb{F}^\times_q=\mathbb{F}_q\setminus \{0\}$ denote the multiplicative groups of $\mathbb{F}_{q^m}$ and $\mathbb{F}_q$, respectively.
Note that the kernel of the trace map from $\mathbb{F}_{q^m}$ to $\mathbb{F}_q$
\begin{equation*}
    K=\left\{x\in \mathbb{F}_{q^m} \mid\operatorname{Tr}_{\mathbb{F}_{q^m}/\mathbb{F}_q}(x)=x+x^q+\cdots+x^{q^{m-1}}=0\right\}
\end{equation*}
is the union of $\{0\}$ and a collection of cosets of $\mathbb{F}^\times_q$ in $\mathbb{F}^\times_{q^m}$.
Consequently, $K\setminus 0$ is a union of full $\mathbb{F}^\times_q$-cosets, and hence we can identify $K\setminus 0$ as a subset of $\mathcal{G}$.
Furthermore, under this identification, $K \setminus 0$ is a difference set of $\mathcal{G}$ with parameters
\begin{equation*}
    (v,k,\lambda)=\left(\frac{q^m-1}{q-1},\frac{q^{m-1}-1}{q-1},\frac{q^{m-2}-1}{q-1} \right),
\end{equation*}
which is defined to be the \emph{Singer difference set} of $\mathcal{G}$.
For $m=2$, this reduces to a trivial $(q+1, 1, 0)$-difference set.

\subsubsection{Paley difference set}
Let $q$ be a prime power with $q\equiv 3 \pmod{4}$. Let $\mathbb{F}_q$ be the finite field with $q$-elements, 
and take $\mathcal{G}$ to be $\mathbb{F}_q^+$, which is the additive group of $\mathbb{F}_q$.
Let 
\begin{equation*}
    D=\{x^2\mid x\in \mathbb{F}_q^\times\}.
\end{equation*}
Then $D$ is a difference set in $\mathcal{G}$ with parameters
\begin{equation*}
    (v,k,\lambda)=(q,\frac{q-1}{2},\frac{q-3}{4}),
\end{equation*}
which is called the \emph{Paley difference set} in $\mathcal{G}$.

\subsection{Basic number theory}
The following lemma is well-known, whose proof we omit.
\begin{lemma}\label{lemma:n1}
    For any positive integer $m,n_1,n_2$ with $m>1$, we have 
    \begin{equation*}
        \operatorname{gcd}(m^{n_1}-1,m^{n_2}-1)=m^{\operatorname{gcd}(n_1,n_2)}-1.
    \end{equation*}
\end{lemma}

The following lemma is a special case of Zsigmondy's theorem \cite{Z92}.
\begin{lemma}\label{lemma:ppd}
    For any integers $m,n$ with $m>1$ and $n>1$, one of the following holds
    \begin{enumerate}[\normalfont(i)]
        \item There is a prime number $\pi$ such that $\pi$ divides $m^n-1$, and $\pi$ does not divide $m^{n'}-1$ for every positive integer $n'<n$;
        \item $n=2$ and $m+1$ is a power of $2$;
        \item $m=2$ and $n=6$.
    \end{enumerate}
\end{lemma}
Furthermore, a prime $\pi$ satisfying condition~\rm(i) in Lemma~\ref{lemma:ppd} is called a \emph{primitive prime divisor} of $m^n-1$.

\begin{lemma}\label{lemma:ppd2}
The following statements hold
\begin{enumerate}[\normalfont (i)]
    \item     For any odd integer $n_1>1$, if $\pi$ is a prime divisor of $2^{n_1}-1$, then $\pi \equiv \pm1 \pmod{8}$.
    \item     Conversely, for every prime $\pi \equiv -1 \pmod{8}$, there exists an odd integer $n_2>1$ such that $\pi$ is a primitive prime divisor of $2^{n_2}-1$.
    Furthermore, for any integer $n_3>1$, $\pi$ divides $2^{n_3}-1$ if and only if $n_2$ divides $n_3$. 
\end{enumerate} 
\end{lemma}
\begin{proof}
    For \rm(i), the condition $\pi$ divides $2^{n_1}-1$ is equivalent to $2^{n_1}=1$ in $\mathbb{F}^\times_{\pi}$, where $\mathbb{F}^\times_{\pi}$ is the multiplicative group of $\mathbb{F}_{\pi}$.
    It follows that $2$ has an odd order dividing $n_1$ in $\mathbb{F}^\times_\pi$. 
    Hence, $2$ is a square in $\mathbb{F}^\times_{\pi}$ and $2$ is a quadratic residue modulo $\pi$, yielding that $\pi \equiv \pm1 \pmod{8}$.
    
    For \rm(ii), we assume that $\pi \equiv -1 \pmod {8}$.
    Then $2$ is a quadratic residue modulo $\pi$, and hence 
    \begin{equation*}
        2^{\frac{\pi-1}{2}} \equiv 1 \pmod{\pi}.
    \end{equation*}
    Let $n_2$ be the order of $2$ in $\mathbb{F}^\times_\pi$.
    Then $n_2$ divides $\frac{\pi-1}{2}$, and hence $n_2$ is odd.
    Furthermore, by the definition following Lemma~\ref{lemma:ppd}, $\pi$ is a primitive prime divisor of $2^{n_2}-1$.
    It remains to prove the last statement. 
    To this end, note that if $n_2$ divides $n_3$ then $2^{n_2}-1$ divides $2^{n_3}-1$ by Lemma~\ref{lemma:n1}.   Conversely, suppose that $n_3$ is an integer such that $\pi$ divides $2^{n_3}-1$.
    Then, by Lemma~\ref{lemma:n1}, $\pi$ divides 
    \begin{equation*}
        \operatorname{gcd}(2^{n_2}-1,2^{n_3}-1)=2^{\operatorname{gcd}(n_2,n_3)}-1.
    \end{equation*}
    On the other hand, $\pi$ is a primitive prime divisor of $2^{n_2}-1$, 
    so $\pi$ does not divide $2^d-1$ for any $d<n_2$, forcing $\operatorname{gcd}(n_2,n_3)=n_2$.
    
    This completes the proof.
\end{proof}

By an argument similar to the proof of Lemma~\ref{lemma:ppd2}, we obtain the following lemma.
\begin{lemma}\label{lemma:ppd3}
    The following statements hold
    \begin{enumerate}[\normalfont (i)]
        \item For any odd integer $n_1>1$, if $\pi$ is an odd prime divisor of $3^{n_1}-1$, then $\pi \equiv \pm1 \pmod{12}$.
        \item Conversely, for every prime $\pi$ with $\pi \equiv -1 \pmod{12}$, there exists an odd integer $n_2>1$ such that $\pi$ is a primitive prime divisor of $3^{n_2}-1$. 
        Furthermore, for any integer $n_3>1$, $\pi$ divides $3^{n_3}-1$ if and only if $n_2$ divides $n_3$.
    \end{enumerate} 
\end{lemma}

\section{New construction of divisible design graphs}\label{sec:construction}

The construction of the distance-regular graphs from Theorem \ref{thm:DRG1} can be modified to produce divisible design graphs.

\begin{construction}\label{constr}
    Let \(G,q,p,S,g,H,r\) be as in the standing assumptions in Section~\ref{subsec:groups}.
    Let $D$ be an $(r,k,\lambda)$-difference set in $N_G(S)/H$. 
    Select $f_1,\cdots, f_k \in N_G(S)$ such that 
    \[D=\{f_1H,\cdots, f_kH \}.\] 
    Let $X$ be the union of double cosets $\bigcup_{i=1}^{k}Hf_igH$.
    Let $\Gamma(q,D)$ be the coset graph $\Gamma(G,H,X)$.
\end{construction}

\begin{remark}\label{rmk}
The construction is independent of the choice of representatives $\{f_1,\cdots,f_k\}$ of the cosets in $D$.
In particular, assume that $f_iH=f_{i'}H$. Then 
\[Hf_i=f_iH=f_{i'}H=Hf_{i'}\]
since $H \mathrel{\unlhd} N_G(S)$, and hence $Hf_igH=Hf_{i'}gH$.

Let $f$ be as in Lemma~\ref{lem:sylownormalizer}. 
Then 
\begin{equation}\label{eq:iso}
N_G(S)/H\cong \langle f \rangle/\langle f^r\rangle \cong \mathbb{Z}/r\mathbb{Z}    
\end{equation}

Under this isomorphism, a set of representatives of $\langle f \rangle/\langle f^r \rangle$ in $\langle f \rangle$ is just a set of representatives of $N_G(S)/H$.

We may select $f_1,\cdots, f_k\in \langle f \rangle$ such that 
\[\{f_1H, \cdots, f_kH\}=D.\]
Furthermore, if $D$ corresponds to
$\{i_1,\cdots,i_k\}\subseteq \mathbb{Z}/r\mathbb{Z}$ under the above isomorphism (\ref{eq:iso}), then we could take $\{f_1, \cdots, f_k\}$ to be $\{f^{i_1},\cdots,f^{i_k}\}$. 

Clearly, $1\notin X$, since $f_ig\notin H$ for every $1\leq i \leq k$. 
By taking $f_1,\cdots,f_k$ from elements in $\langle f \rangle$, we see that $f_ig$ is an involution for all $1\leq i \leq k$ (Lemma \ref{lem:sylownormalizer} \rm(iv)).
It follows that $Hf_igH=(Hf_igH)^{-1}$ and $X=X^{-1}$ in the construction.
\end{remark}

Before we state our main theorem and prove it, we first study some properties of $\Gamma(q,D)$.
\begin{lemma}\label{nbhdH}
    The neighborhood of $H$ is \[\{Hf_igs \mid 1\leq i \leq k, s\in S \},\]
    and $\Gamma(q,D)$ is regular of valency $|S|\cdot k$.
\end{lemma}
\begin{proof}
    By Remark~\ref{rmk}, we may choose the representatives $f_1,\cdots, f_k$ so that $\{f_1,\cdots,f_k\} \subseteq \{1, f,\cdots,f^{r-1} \}$, where $f$ is as in Lemma~\ref{lem:sylownormalizer}.
    By combining the definition of coset graphs which is introduced in Section \ref{sec:cosgraph} and Lemma~\ref{lem:sylownormalizer} \rm(vi), \rm(vii), we have
    \[Hx\sim H \iff x\in \bigcup_{i=1}^{k}Hf_igH=\bigsqcup_{i=1}^kHf_igS=\bigsqcup_{i=1}^k\bigsqcup_{s\in S}Hf_igs.\]
    Thus, $H$ has $|S|\cdot k$ neighbors.
    Finally, $\Gamma(q,D)$ is regular since $G$ acts transitively on its vertices by right multiplication (see Section \ref{sec:cosgraph}).
\end{proof}

The following lemma will be useful in the proof of our main theorem.
\begin{lemma}\label{lemma:intsg}
    For any $0\leq i,j,i' \leq r-1$, let $f^iS^gf^j=\{f^is^gf^j\mid s\in S\}$.
    Then
    \[ |f^iS^gf^j \cap Hf^{i'}gH|=\frac{|S|-1}{r}.\]
\end{lemma}
\begin{proof}
    Note that \[f^iS^gf^j=f^{i+j}(f^{-j}S^gf^j)=f^{i+j}(S^g) \subseteq N_G(S^g)=N_G(S)^g,\]
    where the second equality is by Lemma~\ref{lem:sylownormalizer} \rm(ii).
    Furthermore, $f^is^gf^j\in N_G(S)\cap N_G(S)^g=\langle f \rangle$ if and only if $s=1$, by combining Lemma~\ref{lem:sylownormalizer} \rm(iv) and the fact that $S^g \cap \langle f \rangle=1$. 
    It follows from Lemma~\ref{lem:sylownormalizer} \rm(vi) that 
    \[|f^iS^gf^j \cap N_G(S)|= |f^iS^gf^j \cap \bigcup_{\ell=0}^{r-1}Hf^\ell|=1,\] and hence 
    \[|f^iS^gf^j \cap \bigcup_{\ell=0}^{r-1}Hf^\ell gH|=|f^iS^gf^j|-1=|S|-1\] by Lemma~\ref{lem:sylownormalizer} \rm(vii).
    It suffices to show that for each $0 \leq i',j' \leq r-1$, we have $|f^iS^gf^j \cap Hf^{i'}gH|=|f^iS^gf^j \cap Hf^{j'}gH|.$
    To see this, note that $f\in N_G(S^g)$ by Lemma~\ref{lem:sylownormalizer} \rm(iv), and hence for any integer $t$, we have 
    $(f^iS^gf^j)^{f^t}=(f^i)^{f^t}(S^g)^{f^t}(f^j)^{f^t}=f^iS^gf^j.$ 
    On the other hand, 
    \[(Hf^{i'}gH)^{f^t}=H^{f^t}(f^{i'})^{f^t}g^{f^t}H^{f^t}=Hf^{i'}f^{-t}(gf^tg)gH=Hf^{i'-2t}gH,\]
    because $f$ normalizes $H$ (Lemma~\ref{lem:sylownormalizer} \rm(v)) and $f^g=f^{-1}$ (Lemma~\ref{lem:sylownormalizer} \rm(iv)).
    Furthermore, as $r$ is odd, it follows that there exists an integer $t$ such that $i'-2t\equiv
    j' \pmod{r}$.
    Thus, $f^{i'-2t}\in \langle f^r \rangle f^{j'} \leq Hf^{j'}$, and hence $(Hf^{i'}gH)^{f^t}=Hf^{i'-2t}gH=Hf^{j'}gH$.
    Therefore, 
    \begin{equation*}
\begin{split}
  |f^iS^gf^j\cap Hf^{i'}gH|
  &= |(f^iS^gf^j\cap Hf^{i'}gH)^{f^t}| \\
  &= |(f^iS^gf^j)^{f^t} \cap (Hf^{i'}gH)^{f^t}|
     = |f^iS^gf^j\cap Hf^{j'}gH|,
\end{split}
\end{equation*}
    as required.
\end{proof}

\begin{theorem}\label{thm:DDG}
Let \(G,q,p,S,g,H,r\) be as in the standing assumptions in Section~\ref{subsec:groups}, and let $\Gamma(q,D)$ be as in Construction~\ref{constr}.
Then $\Gamma(q,D)$ is a divisible design graph with parameters
\begin{equation}\label{eq:ddgpar}
((|S|+1)r,|S|\cdot k,|S|\cdot \lambda,\frac{|S|-1}{r}k^2,|S|+1,r).    
\end{equation}
In particular, the following statements hold

\begin{enumerate}[\normalfont(i)]
\item If $G = \Sz(q)$, then the graph $\Gamma(q,D)$ is a divisible design graph with parameters
$$
((q^2+1)r, q^2k, q^2\lambda, \frac{q^2-1}{r}k^2,q^2+1,r).
$$
\item If $G = \Ree(q)$, then the graph $\Gamma(q,D)$ is a divisible design graph with parameters
\[((q^3+1)r, q^3k, q^3\lambda, \frac{q^3-1}{r}k^2,q^3+1,r).\]
\end{enumerate}
\end{theorem}

\begin{proof}
    First, note that by Lemma~\ref{lem:sylownormalizer} \rm(iii)
    \[ |V(\Gamma(q,D))|=|G:H|=|G:N_G(S)|\cdot |N_G(S):H|=(|S|+1)r.\]
    Furthermore, we have $\Gamma(q,D)$ is regular of valency $|S|\cdot k$ by Lemma~\ref{nbhdH}.
    
    We take the canonical partition to be the right cosets of $N_G(S)$ in $G$.
    That is, two vertices $Hx_1$ and $Hx_2$ are in the same class if and only if they belong to the same right coset of $N_G(S)$ in $G$, which is further equivalent to $x_1x_2^{-1}\in N_G(S)$.
    Note that $|N_G(S):H|=r$, it follows from Lemma~\ref{lem:sylownormalizer} \rm(iii) that our canonical partition has $|S|+1$ classes, with each class of size $r$.

    Recall the transitive action of $G$ on the points of $\Gamma(q,D)$ by right multiplication (Lemma~\ref{lem:trans}).
    We next claim that
    
    \noindent{\emph{Claim: $G$ acts doubly transitively by right multiplication on the classes of our canonical partition}}

    To see this, note that the stabilizer of the class containing $H$ is $N_G(S)$.
    Since $N_G(S)$ is also the stabilizer of the Sylow $p$-subgroup $S$ under the conjugation action of $G$, and $G$ acts doubly transitively on its Sylow $p$-subgroups by conjugation, the action of $G$ on the classes is permutation equivalent to that doubly-transitive action by \cite[Lemma 1.6B]{PGD96}. 
    Hence, $G$ acts doubly transitively on the classes.
    This completes our claim.
    
    Alternatively, from (\ref{eqngs}) in the proof of Lemma~\ref{lem:sylownormalizer}, one sees that the nontrivial right cosets of $N_G(S)$ in $G$ are exactly $N_G(S)gs$ with $s\in S$.
    It follows that $S$ acts transitively on them by right multiplication.
    Since $S\leq N_G(S)$, this again shows that $G$ acts doubly transitively on the classes by right multiplication.
     
    We then proceed to determine the number of common neighbours of two vertices $v_1,v_2\in \Gamma(q,D)$, we deal with the cases where $v_1,v_2$ are in the same class and $v_1,v_2$ are in different classes separately.

    \noindent{\emph{Case 1: $v_1,v_2$ are in the same class}}
    
    Note that $G$ acts transitively on $V(\Gamma(q,D))$ by right multiplication (Lemma~\ref{lem:trans}).
    We may replace $v_1$ and $v_2$ with $v_1g$ and $v_2g$ for a suitable $g$ and assume that $v_1=H$.  
    Then $v_2=Hf^j$ for some $1\leq j \leq r-1$, since $v_1$ and $v_2$ lie in the same class, i.e. the same right coset of $N_G(S)$ in $G$, and $N_G(S)=\bigsqcup_{j=0}^{r-1}Hf^j$.

    Let $v_3\in V(\Gamma(q,D))$ be adjacent to $H$. Then $v_3=Hf_igs$ for some $1\leq i \leq k$ and some $s\in S$ by Lemma~\ref{nbhdH}. 
    By the definition of $\Gamma(q,D)$ as a coset graph, we have $v_3\sim v_2$ if and only if $f_igsf^{-j}\in \bigcup_{i'=1}^{k}Hf_{i'}gH$. On the other hand,
   \[f_igsf^{-j}=f_i(gf^{-j}g)gs^{f^{-j}}=f_if^jgs^{f^{-j}}\in Hf_if^jgH,\]
    where $gf^{-j}g=f^j$ by Lemma~\ref{lem:sylownormalizer} \rm(iv) and $s^{f^{-j}}\in S\leq H$.
    It follows from Lemma~\ref{lem:sylownormalizer} \rm(vii) that $v_3$ is adjacent to $v_2$ if and only if $f_if^j\in Hf_{i'}$ for some $1\leq i'\leq k.$
    Since $\{f_1H,\cdots,f_kH\}=D$ is an $(r,k,\lambda)$-difference set of $N_G(S)/H \cong \mathbb{Z}/r\mathbb{Z}$, 
    it follows that there are exactly $\lambda$ possible $i\in [k]$ such that $f_if^j\in \bigcup_{i'=1}^kHf_{i'}$. 
    As a result, among the neighbors $\{Hf_igs \mid 1\leq i\leq k,s\in S\}$ of $v_1=H$, there are exactly $|S| \cdot \lambda$ vertices that are also adjacent to $v_2$. 
    
    This completes Case 1.

    \noindent{\emph{Case 2: $v_1$ and $v_2$ are in different classes}}

    By the double transitivity established in the claim, we may replace $v_1$ and $v_2$ with $v_1g_1$ and $v_2g_1$ for some suitable $g_1\in G$ and assume that $v_1$ lies in the class $N_G(S)$ and $v_2$ lies in the class $N_G(S)g$.
    Moreover, the class containing $H$ consists of vertices $Hx$ with $x\in N_G(S)$, and the class containing $Hg$ consists of vertices $Hx$ with $x\in N_G(S)g$. 
    By a direct computation and Lemma~\ref{lem:sylownormalizer} \rm(iv), the stabilizer of both classes in $G$ acting by right multiplication is $N_G(S)\cap N_G(S)^g=\langle f \rangle$. 
    Since $\langle f \rangle$ acts transitively on the vertices in the class containing $H$ (which are just $\{H, Hf,\cdots, Hf^{r-1}\}$) by right multiplication, we may further replace $v_1$ and $v_2$ with $v_1f^t$ and $v_2f^t$ for some integer $t$ and thereby assume that $v_1=H$.

    Thus, by our assumption in the last paragraph, we have $v_1=H$, and $v_2=Hf^jg$ for some integer $0\leq j \leq r-1$. 
    For a vertex $v_3$ in $\Gamma(q,D)$ adjacent to $v_1$, we have $v_3=Hf_{i'}gs$ for some $1\leq i'\leq k$ and some $s\in S$ by Lemma~\ref{nbhdH}.
    In view of Remark~\ref{rmk}, we may choose the representatives $\{f_1,\cdots,f_k\}\subseteq \{1,f,\cdots,f^{r-1}\}$.
    By the construction of $\Gamma(q,D)$, we have that $v_3$ is adjacent to $v_2$ if and only if 
    \[f_{i'}gsgf^{-j}\in \bigcup_{i=1}^{k}Hf_igH=\bigsqcup_{i=1}^{k}Hf_igH,\] 
    where the last equality holds by Lemma~\ref{lem:sylownormalizer} \rm(vii).
    Furthermore, $f_{i'}gsgf^{-j}=f_{i'}s^gf^{-j}$ with $s^g\in S^g$ and $f_{i'}\in \{1,f,\cdots,f^{r-1}\}$ by our choice.
    It follows from Lemma~\ref{lemma:intsg} that there are exactly $(|S|-1)k/r$ elements $s\in S$ such that $f_{i'}gsgf^{-j}\in \bigcup_{i=1}^{k} Hf_igH$ for fixed $i'$. 
    Therefore, among all the neighbors $\{Hf_igs \mid 1\leq i \leq k, s\in S\}$ of $H$ (Lemma~\ref{nbhdH}), there are exactly $(|S|-1)k^2/r$ of them that are also adjacent to $v_2$ as $i'$ runs from $1$ to $k$.
    This completes Case 2.

    We have already calculated all the parameters of $\Gamma(q,D)$ in (\ref{eq:ddgpar}) as a divisible design graph.
    Finally, we substitute the value of $|S|$ (see Lemma~\ref{lem:sylownormalizer} \rm(i)) into (\ref{eq:ddgpar}) and obtain \rm(i) and \rm(ii).

    This completes the proof.
\end{proof}

Finally, we present examples for Construction~\ref{constr} and Theorem~\ref{thm:DDG}.

\begin{example}
Take $D$ as the trivial difference set $\{H\}$, where $H$ denotes the identity element of the cyclic group $N_G(S)/H$. 
Then the graph given in Construction~\ref{constr} is just the graph in Tsiovkina's construction (Theorem~\ref{thm:DRG1}).
\end{example}
The above example shows that Tsiovkina's construction is a special case of our construction.
Next, we present some examples by taking $D$ to be non-trivial cyclic difference sets in Construction~\ref{constr}. We start by examples of plugging Singer difference sets.

\begin{example}
    Let $G=\Sz(q)$, where $q=2^{2n+1}$ for some positive integer $n$. 
    Let $n_1, n_2$ be divisors of $2n+1$ such that $n_2\geq 2n_1\geq 2$ and $n_1$ divides $n_2$.
    Then $2^{n_1}-1$ divides $2^{n_2}-1$, and $2^{n_2}-1$ divides $2^{2n+1}-1$ by Lemma~\ref{lemma:n1}.
    Let $H$ be the unique subgroup of index $r=\frac{2^{n_2}-1}{2^{n_1}-1}$ in $N_G(S)$ (see Lemma~\ref{lem:sylownormalizer} \rm(v) for the structure of $H$).
    Take $D$ to be the Singer difference set of parameters 
    \[(r,k,\lambda)=\left(\frac{2^{n_2}-1}{2^{n_1}-1}, \frac{2^{n_2-n_1}-1}{2^{n_1}-1}, \frac{2^{n_2-2n_1}-1}{2^{n_1}-1}\right)\] in the cyclic group $N_G(S)/H \cong \mathbb{Z}/r\mathbb{Z}$.
    Then the graph in Construction~\ref{constr} is a divisible design graph of parameters
    \begin{equation*}
\begin{aligned}
&\left(
\frac{(2^{4n+2}+1)(2^{n_2}-1)}{2^{n_1}-1},\,
\frac{2^{4n+2}(2^{n_2-n_1}-1)}{2^{n_1}-1},\,
\frac{2^{4n+2}(2^{n_2-2n_1}-1)}{2^{n_1}-1},
\right.\\
&\left.
\frac{(2^{4n+2}-1)(2^{n_2-n_1}-1)^2}{(2^{n_1}-1)(2^{n_2}-1)},\,
2^{4n+2}+1,\,
\frac{2^{n_2}-1}{2^{n_1}-1}
\right).
\end{aligned}
\end{equation*}
\end{example}

\begin{example}
    Let $G=\Ree(q)$, where $q=3^{2n+1}$ for some positive integer $n$. 
    Let $n_2 \geq 2n_1 \geq 2$ such that $n_1$ divides $n_2$ and $n_2$ divides $2n+1$.
    Then $3^{n_1}-1$ divides $3^{n_2}-1$ and $3^{n_2}-1$ divides $3^{2n+1}-1$ by Lemma~\ref{lemma:n1}.
    Let $H$ be the unique subgroup of index $r=\frac{3^{n_2}-1}{3^{n_1}-1}$ in $N_G(S)$, where 
    \[r=1+3^{n_1}+\cdots +3^{(\frac{n_2}{n_1}-1)n_1}\] is odd.
    We take $D$ to be the Singer difference set of parameters 
    \begin{equation*}
        (r,k,\lambda)=\left(\frac{3^{n_2}-1}{3^{n_1}-1},\frac{3^{n_2-n_1}-1}{3^{n_1}-1},\frac{3^{n_2-2n_1}-1}{3^{n_1}-1} \right)
    \end{equation*}
    in the cyclic group $N_G(S)/H\cong \mathbb{Z}/ r\mathbb{Z}$.
    Then the graph in Construction~\ref{constr} is a divisible design graph of parameters
    \begin{equation*}
\begin{aligned}
&\left(\frac{(3^{6n+3}+1)(3^{n_2}-1)}{3^{n_1}-1},\,
\frac{3^{6n+3}(3^{n_2-n_1}-1)}{3^{n_1}-1},\,
\frac{3^{6n+3}(3^{n_2-2n_1}-1)}{3^{n_1}-1},\right.\\
&\left.\frac{(3^{6n+3}-1)(3^{n_2-n_1}-1)^2}{(3^{n_1}-1)(3^{n_2}-1)},\,
3^{6n+3}+1,\, \frac{3^{n_2}-1}{3^{n_1}-1}\right).
\end{aligned}
\end{equation*}
     
\end{example}

Next, we present examples of new divisible design graphs by plugging Paley difference sets into Tsiovkina's construction.
Note that in our construction (Construction~\ref{constr}), $D$ can only be a cyclic difference set. 
On the other hand, a Paley difference set is cyclic if and only if the group is the additive group of $\mathbb{F}_{\pi}$, where $\pi \equiv 3 \pmod 4$ is a prime.
By Lemma~\ref{lemma:ppd2} (resp. Lemma~\ref{lemma:ppd3}), if a prime $\pi \equiv 3 \pmod{4}$ divides $2^{2n+1}-1$ (resp. $3^{2n+1}-1$) with $n\geq 1$, then $\pi\equiv -1 \pmod{8}$ (resp. $\pi\equiv -1 \pmod{12}$).
We next construct examples by plugging the Paley difference set of the additive group of $\mathbb{F}_{\pi}$ into Construction~\ref{constr}, for primes  $\pi \equiv -1 \pmod{8}$ (Suzuki case) and primes $\pi \equiv -1 \pmod{12}$ (Ree case).

\begin{example}
    Let $\pi \equiv -1 \pmod{8}$ be a prime, and let $n_2$ be an odd integer such that $\pi$ is a primitive prime divisor of $2^{n_2}-1$. Such an integer exists by Lemma~\ref{lemma:ppd2}. 
    By Lemma~\ref{lemma:ppd2} \rm(ii), we have $\pi$ divides $2^{2n+1}-1$ if and only if $n_2$ divides $2n+1$. 
    We then take $n$ as an integer such that $n_2$ divides $2n+1$.
    Finally, let $G=\Sz(q)$ with $q=2^{2n+1}$, and let $H$ be the unique subgroup of index $r=\pi$ in $N_G(S)$.
    Take $D$ to be the Paley difference set in $N_G(S)/H\cong \mathbb{Z}/\pi \mathbb{Z}$ with parameters
    \begin{equation*}
        (r,k,\lambda)=\left(\pi, \frac{\pi-1}{2}, \frac{\pi-3}{4}\right).
    \end{equation*}
    Then the graph in Construction~\ref{constr} is a divisible design graph of parameters
    \begin{equation*}
        \left((2^{4n+2}+1)\pi, \frac{2^{4n+2}(\pi-1)}{2}, \frac{2^{4n+2}(\pi-3)}{4}, \frac{(2^{4n+2}-1)(\pi-1)^2}{4\pi},2^{4n+2}+1,\pi \right).
    \end{equation*}
\end{example}

\begin{example}
    Let $\pi \equiv -1 \pmod{12}$ be a prime, and let $n_2$ be an odd integer such that $\pi$ is a primitive prime divisor of $3^{n_2}-1$. 
    Such an integer exists by Lemma~\ref{lemma:ppd3}. 
    By Lemma~\ref{lemma:ppd3} \rm(ii), we have $\pi$ divides $3^{2n+1}-1$ if and only if $n_2$ divides $2n+1$. 
    We then take $n$ as an integer such that $n_2$ divides $2n+1$.
    Finally, let $G=\Ree(q)$ with $q=3^{2n+1}$, and let $H$ be the unique subgroup of index $r=\pi$ in $N_G(S)$.
    Take $D$ to be the Paley difference set in $N_G(S)/H\cong \mathbb{Z}/\pi \mathbb{Z}$ with parameters
    \begin{equation*}
        (r,k,\lambda)=\left( \pi, \frac{\pi-1}{2}, \frac{\pi-3}{4} \right).
    \end{equation*}
    Then the graph in Construction~\ref{constr} is a divisible design graph of parameters
    \begin{equation*}
        \left((3^{6n+3}+1)\pi, \frac{3^{6n+3}(\pi-1)}{2}, \frac{3^{6n+3}(\pi-3)}{4}, \frac{(3^{6n+3}-1)(\pi-1)^2}{4\pi},3^{6n+3}+1,\pi \right).
    \end{equation*}
\end{example}

In summary, by plugging Singer difference sets and Paley difference sets into Tsiovkina's construction, we obtain four infinite families of new divisible design graphs.
\section*{Acknowledgments}

\noindent Sergey Goryainov is supported by the Natural Science Foundation of Hebei Province (A2023205045) and the 111 Center (Grant No.D26018). Yilin Xie is supported by Science Foundation of Hebei Normal University (L2026B48).

\end{document}